\documentclass[12pt]{article}
\usepackage{amssymb,latexsym, amsmath, amsthm, amsfonts, amssymb, enumerate, paralist,cite}
\usepackage[notref,notcite]{showkeys}
\usepackage[usenames,dvipsnames,svgnames,x11names, table]{xcolor}
\usepackage{graphicx, overpic}
\usepackage{color}

\usepackage{hyperref}
\hypersetup{
                CJKbookmarks=true
	colorlinks=true,
	linkcolor=blue,
	filecolor=blue,
	urlcolor=blue,
	citecolor=blue,
}
\def\H{\mathbb H}
\def\UBB{\mathbb U}
\def\R{\mathbb R}
\def\C{\mathbb C}
\def\N{\mathbb N}
\def\Z{\mathbb Z}
\def\X{{\mathrm{X}}}
\def\xx{\mathcal{X}}
\def\G{\mathbb G}
\def\M{\mathbb M}

\def\Cq{(2\pi)^\frac{q}{2}}
\def\tCqm{\widetilde{C}(q,m)}
\def\Cqm{C(q,m)}
\def\OA{\Omega_*}
\def\OM{\widetilde{\M}}
\def\rPo{\mathrm{P}_h^{\{0\}}}
\def\rPk{\mathrm{P}_h^{\{k\}}}
\def\rPkk{\mathrm{P}_h^{\{k+1\}}}
\def\rQ{\mathrm{Q}}
\def\rR{\mathrm{R}}
\def\rF{\mathfrak{F}}
\def\cZ{\mathcal{Z}}
\def\OAk{\Omega_*^{(k)}}
\def\OAko{\Omega_*^{(k_0)}}
\def\OAkp{\Omega_*^{(k + 1)}}
\def\rFks{\phi_k}
\def\rFzs{\phi_0}
\def\rFkps{\phi_{k + 1}}
\def\srFks{\phi_k^*}
\def\srFkps{\phi_{k + 1}^*}
\def\srFzs{\phi_0^*}
\def\isrFzs{\Phi_0}
\def\isrFks{\Phi_k}
\def\isrFkzs{\Phi_{k_0}}
\def\isrFkps{\Phi_{k + 1}}
\def\rFk{\phi^{\{k\}}}
\def\Mk{\mathrm{M}^{(k)}}
\def\ss{\mathbf{s}}
\def\He{\mathrm{Hess}}

\newtheorem{theo}{Theorem}[section]
\newtheorem{lem}{Lemma}[section]

\newtheorem{prop}{Proposition}[section]

\newtheorem{alg}{Algorithm}[section]

\title{A complete answer to the Gaveau--Brockett problem}
\author{Hong-Quan Li, Ye Zhang}
\date{}
\begin{document}

\renewcommand{\theequation}{\arabic{equation}}
\setcounter{equation}{0} \maketitle

\vspace{-1.0cm}

\bigskip

{\bf Abstract.}
The note is dedicated to provide a satisfying and complete answer to the long-standing Gaveau--Brockett open problem. More precisely, we determine the exact formula of the Carnot--Carath\'eodory distance on arbitrary step-two groups. The basic idea of the proof is combining Varadhan's formulas with the explicit expression for the associated heat kernel $p_h$ and the method of stationary phase. However, we have to introduce a number of original new methods, especially the usage of the concept of ``Operator convexity''. Next, new integral expressions for $p_h$ by means of properties of Bessel functions will be presented. An unexpected direct proof for the well-known positivity of $p_h$ via its original integral formula, will play an important role. Furthermore, all normal geodesics joining the identity element $o$ to any given $g$  as well as the cut locus can be characterized on every step-two groups. Finally, the corresponding results in Riemannian geometry on step-two groups will be briefly presented as well.

\medskip


\medskip

{\bf Significance. }
The most fundamental problem in geometry is to determine the distance as well as the shortest path(s) between any two given points. Another basic problem is to find the cut locus. To solve these problems is impractical in general. In the setting of step-two groups, the first problem in the sub-Riemannian context is well-known as the Gaveau--Brockett problem. Despite a lot of related works, a complete description for these basic geometric objects is obtained in few cases. Surprisingly, we can solve the Gaveau--Brockett problem by a fully analytic approach, specifically via combining tools from Convex, Fourier, Matrix Analysis as well as the heat kernel. In addition, our method can be adapted to study the asymptotic behavior at infinity of the heat kernel.

\medskip

\section*{Introduction}

In the past several decades, step-two groups and their sub-Laplacians, as special Lie groups of polynomial volume growth or perfect  sub-Riemannian manifolds, have attracted wide attention from experts in various fields, such as complex analysis, control theory, geometric measure theory, harmonic analysis, heat kernel, Lie group theory, probability analysis, PDE, sub-Riemannian geometry, etc. The literature on these topics is huge.  We only mention several closely related monographs here: \cite{B84, VSC92, M02, BLU07, R14, ABB20}.

In this note, we will mainly focus on the sub-Riemannian geometry of 2-step groups. Despite recent progresses on sub-Riemannian geometry,   some most fundamental problems are far from being solved, or even poorly known, in this very fine framework. From the left-invariant property, these problems include determining the Carnot--Carath\'eodory (or sub-Riemannian) distance $d(g) := d(g, o)$ between $g$ and the identity $o$  (in other words, the Gaveau--Brockett problem, cf. \cite{G77}, \cite{B82}, and \cite{Br84}), and the cut locus  $\mathrm{Cut}_o$  of $o$ (namely, the set of points where the squared distance is not smooth).

It has been known for a long time that behaviors of the heat kernel are closely related to geometric properties of the underlying space, see \cite{V67}, \cite{M75}, \cite{A81}, and \cite{B84} for more information about early works. In particular, the  Varadhan's formula holds in our setting of step-two groups.  In addition, an explicit expression for the heat kernel, by means of oscillatory integral, is well-known (cf.  \eqref{2c0} below). The standard method to treat the asymptotic behavior of the oscillatory integral is that of stationary phase.  Naturally, we study the distance by combining the Varadhan's formula with the method  of stationary phase. However, this task is by no means easy, and  a number of original new methods in this context  are indispensable.

In this work we omit details of the proofs and prefer to insist on the main stream of ideas. We plan to include the detailed proofs  as well as the main results of \cite{Li20, LZ21} in a single work in the future.

\section{Preliminaries} \label{ss1}

Recall that a step-two Carnot group (or step-two group, 2-step group, in short) can be considered as $\R^q \times \R^m$, $q, m \in \N^* = \{1, 2, 3, \ldots\}$, with the group law
\[
(x , t) \cdot (x^{\prime}, t^{\prime}) = \left(x + x^{\prime}, t + t^{\prime} + \frac{1}{2}\langle  \UBB \, x, x^{\prime} \rangle \right), \quad
g := (x, t) \in \R^q \times \R^m,
\]
where $\langle\UBB \, x,  x^{\prime} \rangle := (\langle U^{(1)} x, x^{\prime} \rangle, \ldots, \langle   U^{(m)} x,x^{\prime} \rangle) \in \R^m$. Here $\UBB = \{U^{(1)},\ldots,U^{(m)}\}$ is an $m$-tuple of linearly independent $q \times q$ skew-symmetric matrices with real entries and $\langle \cdot, \cdot \rangle$ (or $\mbox{} \cdot \mbox{}$ in the sequel when there is no ambiguity) denotes the usual inner product on $\R^q$.  The Haar measure $dg$ coincides with the $(q + m)$-dimensional Lebesgue measure on $\G$. Let $U^{(j)} = (U^{(j)}_{l, k})_{1 \leq l, k \leq q}$ ($1 \leq j \leq m$). The canonical sub-Laplacian is given by
\begin{align*}
\Delta := \sum_{l = 1}^q \X_l^2, \quad \mbox{with} \quad \X_l(g) : = \frac{\partial}{\partial x_l} + \frac{1}{2} \sum_{j = 1}^m \left( \sum_{k = 1}^{q} U^{(j)}_{l, k} \, x_k \right) \frac{\partial}{\partial t_j}, \quad 1 \leq l \leq q.
\end{align*}
It is well-known that the heat semigroup $(e^{h \Delta})_{h > 0}$ has a convolution kernel $p_h$, in the sense that
\begin{eqnarray*}
e^{h \Delta} u(g) = u \ast p_h(g) = \int_{\G} u(g') \, p_h((g')^{-1} \cdot g) \, dg' \quad \mbox{for suitable functions $u$}.
\end{eqnarray*}
The following explicit formulas can be found in \cite{C79} or \cite[Theorem 4, p.\,293]{BGG96} with a slight modification of notations:
\begin{align} \label{2c0}
p_h(x, t) = \frac{\Cqm}{h^{\frac{q}{2} + m}} \, \int_{\R^m} V(\lambda) \,  \exp\left\{ -\frac{\widetilde{\phi}((x, t); \lambda)}{4h} \right\} \, d\lambda, \qquad \forall \, h > 0, \, (x,t) \in \G,
\end{align}
where the  positive constant $\Cqm$ depends only on $q$ and $m$, \begin{gather*}
V(\lambda) := \left( \mathrm{det} \frac{U(\lambda)}{\sinh{U(\lambda)}} \right)^{\frac{1}{2}},  \quad
\widetilde{\phi}((x, t); \lambda) := \langle U(\lambda) \coth{U(\lambda)} \, x, \, x  \rangle - 4  i \, t \cdot \lambda,
\end{gather*}
with $ U(\lambda) := i \, \widetilde{U}(\lambda) := i  \sum\limits_{j = 1}^m \lambda_j \, U^{(j)}$.

In the sequel, $\|  U(\lambda)  \|$ denotes the usual norm of the linear operator $U(\lambda)$ on $\C^q$, $|x| = \sqrt{x \cdot  x}$, $0$ the number $0$ or the origin in Euclidean spaces, and $\nabla_\theta = \left(\frac{\partial}{\partial \theta_1}, \ldots, \frac{\partial}{\partial \theta_m}\right)$  the usual gradient on $\R^m$. Moreover, we will use $\exp: T_o^* \G \cong \R^q \times \R^m \to \G$ to denote the sub-Riemannian exponential map based at the identity $o$. In other words, $\exp{(\zeta, \tau)}$ is the extremity at time $1$ of the normal geodesic starting from $o$ with the initial covector $(\zeta, \tau)$, which is defined as the projection on $\G$ of the corresponding integral curve of the sub-Riemannian Hamiltonian on the cotangent bundle. See for example \cite[\S~2.3]{R14} and \cite[\S~13.1]{ABB20}  for more details.
Recall that the Varadhan's formula is valid in our situation (cf. \cite{L87}, \cite{BA88}, \cite{VSC92}, and the references therein):
\begin{align} \label{VF}
\lim_{h \longrightarrow 0^+} 4 h \ln{p_h(g)} = -d(g)^2, \qquad \forall \, g \in \G.
\end{align}

\section{Main results and outlines of the proof}

\subsection{An enlightening example and some general results} \label{ss2}

Before going further, we present  the simplest example which will help explain what we have in mind. Let $\H_3$ denote the Heisenberg group, namely $q = 2$, $m = 1$, and $\UBB$ has only one element, which is the matrix of the standard symplectic bilinear form on $\R^2$. Let us introduce some notations in the context of step-two groups, which become very concise in the special case $\H_3$ and will be enclosed in parentheses and marked in blue for convenience. We define the initial reference set as the absorbing, convex, bounded open set
\[
\OA := \left\{ \tau \in \R^m; \ \| U(\tau) \| < \pi \right\} {\color{blue} (= (-\pi, \, \pi))},
\]
and the reference function based at $g = (x,t)$ by
\[
\phi(g; \tau) := \widetilde{\phi}(g; i \tau) = \langle U(\tau) \cot{U(\tau)} \, x, \, x  \rangle + 4 \, t \cdot \tau \ {\color{blue} (=  \tau \cot{\tau} \, |x|^2 + 4 t \tau )}, \qquad \tau \in \OA.
\]

Notice that the function $\phi(g; \cdot)$ is well-defined provided the spectrum of $U(\tau)$ does not contain any $k \, \pi$ ($k \in \Z \setminus \{ 0 \}$), namely $\tau \in  \mathcal{V}$ with
\begin{align} \label{nD1V}
\mathcal{V} := \left\{ \vartheta \in \R^m; \,  \mathrm{det} (k \, \pi \, \mathbb{I}_q - U(\vartheta)) \neq 0, \quad \forall \, k \in \N^* \right\} {\color{blue}(= \R \setminus \{k \pi; \, k \in \Z \setminus \{ 0 \} \}) }.
\end{align}

It follows from \cite[Lemme 3, p. 112]{G77} that the function 
\[
\mu(s) := -\frac{d}{ds} (s \cot{s}) = \frac{2 s - \sin{(2 s)}}{2 \sin^2{s}}
\] 
is an odd function, and a monotonely increasing diffeomorphism between $(-\pi, \, \pi)$ and $\R$. Hence, our first observation on $\H_3$ is the following:

(P1) For every $g$, the reference function $\phi(g; \cdot)$ is smooth and concave on $\OA$. Hence,  $\theta$ is a critical point of $\phi(g; \cdot)$ in $\OA$ if and only if it is a global maximizer in $\OA$.

Furthermore, the following properties on $\H_3$ can be found or deduced directly from \cite{G77, BGG00}:

(P2) Up to a subset of measure zero
\[
\mathcal{W} := \left\{\exp(\zeta_0,  2 \, \theta_0); \ \zeta_0 \in \R^q, \ \theta_0 \in \mathcal{V}^c \right\} \  {\color{blue} (= \{(x, t); \, x = 0 \})},
\]
we can find all the geodesics joining $o$ to any given $g_0  \in \mathcal{W}^c$. To be more precise,
$\exp{(\zeta_0, 2 \, \theta_0) } = g_0 = (x_0, t_0)$ with $\theta_0 \in \mathcal{V}$ if and only if $\nabla_{\theta_0} \phi(g_0; \theta_0) = 0$, i.e.  
\[
t_0 = - \frac{1}{4} \nabla_{\theta_0} \langle U(\theta_0) \, \cot{U(\theta_0)} \, x_0, \ x_0 \rangle \   {\color{blue} (4 t_0 = \mu(\theta_0) |x_0|^2,  \ \theta_0 \not\in \{k \pi; \, k \in \Z \setminus \{ 0 \} \}),}
\]
and
\begin{align*}
x_0 = \left( \frac{U(\theta_0)}{\sin{U(\theta_0)}} e^{-\widetilde{U}(\theta_0)} \right)^{-1} \zeta_0  \quad
{\color{blue}   \mbox{\bigg($x_0 = \frac{\sin{\theta_0}}{\theta_0} \left(
                                                               \begin{array}{cc}
                                                                 \cos{\theta_0} &  - \sin{\theta_0} \\
                                                                 \sin{\theta_0} & \cos{\theta_0} \\
                                                               \end{array}
                                                             \right) \zeta_0$ \bigg).}}
\end{align*}

(P3) Let
\begin{align*}
\M &:=  \left\{ g = (x, t); \exists \, \theta \in \OA \mbox{ s.t. $\theta$ is a nondegenerate critical point (maximizer) of $\phi(g; \cdot)$ in $\OA$} \right\}  \\
&= \left\{  g = \big( x, -\frac{1}{4} \nabla_{\theta} \langle U(\theta) \cot{U(\theta)} \, x, \, x  \rangle \big);  \, x \neq 0, \ \theta \in \OA, \  \det \left(- \mathrm{Hess}_\theta  \, \phi(g; \theta) \right) > 0 \right\} \\
&{\color{blue}  ( = \{(x, t); \, x \neq 0\}) }, \\
\widetilde{\M} &:= \left\{  g = \big( x, -\frac{1}{4} \nabla_{\theta} \langle U(\theta) \cot{U(\theta)} \, x, \, x  \rangle \big);  \, x \in \R^q, \ \theta \in \OA \right\} {\color{blue}( = \H_3 \setminus \{(0, t); \, t \neq 0\})}, \\
\OM_1 &:= \{g; \, \mbox{the set of  global maximizers of $\phi(g; \cdot)$ in $\OA$ is a singleton} \} {\color{blue}( = \{(x, t); \, x \neq 0\})}, \\
\OM_2 &:= \{g; \, \mbox{the set of  global maximizers of $\phi(g; \cdot)$ in $\OA$ has at least two points} \} {\color{blue}( = \{ o \})}.
\end{align*}
It holds that $\overline{\M} = \overline{\widetilde{\M}}$, $\M = \OM_1 \subset \mathcal{S} := \mathrm{Cut}_o^c =  \{g; \, \mbox{$d^2$ is $C^{\infty}$ in a neighborhood of $g$}\}$, and $\OM_2 \subset \mathrm{Cut}_o$.

(P4) $d(g)^2  \geq \sup\limits_{\tau \in \OA} \phi(g; \tau)$  for  every $g \in \G$.

(P5) Let $g = \big( x, -\frac{1}{4} \nabla_{\theta} \langle U(\theta) \cot{U(\theta)} \, x, \, x  \rangle \big) \in \M$ with $\theta \in \OA$. Then the following small-time asymptotics of the heat kernel
\begin{align}\label{smalltM0}
p_h(g) = \Cqm \,  \frac{(8 \pi)^{\frac{m}{2}}}{h^{\frac{q + m}{2}}} \, V(i \theta) \, e^{-\frac{\phi(g; \theta)}{4 h}} \left[  \det \left(- \mathrm{Hess}_\theta  \, \phi(g; \theta) \right) \right]^{-\frac{1}{2}} \left( 1 + O(h) \right), \qquad h \to 0^+,
\end{align}
holds uniformly on every compact set in $\M$.

(P6) If the reference function $\phi(g; \cdot)$ attains its maximum  in $\OA$ at some point, saying $\theta$,
then $d(g)^2 = \phi(g; \theta)$. In other terms,  if $g = (x, t) = \big( x, -\frac{1}{4} \nabla_{\theta} \langle U(\theta) \cot{U(\theta)} \, x, \, x  \rangle \big)$ {\color{blue} ($= (x, \frac{1}{4} \mu(\theta) |x|^2)$)} with $\theta \in \OA$, then we have
\begin{align} \label{DEP}
d(g)^2 = \phi(g; \theta) = \left| \frac{U(\theta)}{\sin{U(\theta)}} x \right|^2 = \max_{\tau \in \OA} \phi(g; \tau) \, {\color{blue} \left( = \big(\frac{\theta}{\sin{\theta}}\big)^2 |x|^2\right)}.
\end{align}

(P7) We have $d(g)^2 = \sup\limits_{\tau \in \OA} \phi(g; \tau)$ provided $g \in \overline{\M}$.

(SP) On $\H_3$, we have that $d(g)^2  = \sup\limits_{\tau \in \OA} \phi(g; \tau)$  for  every $g \in \H_3$,  $\mathrm{Cut}_o = \M^c$,  and the classical cut locus of $o$  (namely, the set of points where geodesics starting at $o$ cease to be shortest)  $\mathrm{Cut}^{\mathrm{CL}}_o = \OM^c$.

We now extend these observations to the general setting:

\begin{theo}
All these properties (P1)-(P7) above are in fact valid on every step-two group.
\end{theo}

More precisely, the claim $\overline{\M} = \overline{\widetilde{\M}}$ in (P3) is exactly \cite[Proposition 4]{LZ21}, whose proof is based on the fact that the function $(x, \, \theta) \mapsto \det(-\mathrm{Hess}_\theta \, \langle U(\theta) \, \cot{U(\theta)} \, x, \ x\rangle)$ is real analytic in $\R^q \times \OA$.  Other results can be
detected in \cite{Li20}. In particular, to get the property (P7) from (P6) or \eqref{DEP} with $g \in \M$, we need the following result:

\begin{prop} \label{NPA}
Let $\Omega \subset \R^m$ be an open, bounded, convex set. Assume that the function $h: \, \R^k \times \Omega \longrightarrow \R$
is continuous and for any $\xx \in \R^k$, $h(\xx, \cdot)$ is concave in $\Omega$. Then the function defined by $H(\xx) := \sup\limits_{\tau \in \Omega} h(\xx; \tau)$
is continuous on $\R^k$.
\end{prop}

To show the remaining properties,  we introduce the usage of  the concept of ``operator convexity''  in this context. To be more specific, the fact that the even function $-s \cot{s}$ is operator convex on $(-\pi, \, \pi)$ plays a crucial role. Additionally, we make use of: (i) the Varadhan's formulas in the proof of (P4); (ii) the method of stationary phase in the proof of (P5); (iii) for (P2), the Spectral Theorem in $\C^q$ as well as the explicit expression of the sub-Riemannian exponential map $\exp$ in the setting of step-two groups, namely the equation of normal geodesics that begin at the identity, which can be found in many works, e.g. \cite{BGG96} with a slight modification.

It is worthwhile to point out that (P6) can be deduced by combining (P1) with (P2) and (P4). Furthermore, under the additional condition that $g \in \M$, \eqref{DEP} is a direct consequence of  \eqref{smalltM0}
and the Varadhan's formulas.

Further results and applications are provided in \cite{Li20} and \cite{LZ21}.
Especially, inspired by (SP), we introduce:

\subsection{GM-groups and their equivalent characterizations}

A step-two group $\G$ is said to be a GM-group (or of type GM) if its subset $\M$, defined in (P3), is dense in $\G$. GM-groups form an enormous class of step-two groups, including many important groups, such as generalized Heisenberg-type groups, step-two groups of corank $1$ or $2$, of Kolmogorov type, and those associated to quadratic CR manifolds. Indeed, all these examples can be explained by the algebraic sufficient condition \cite[Proposition 2.2]{Li20}. Furthermore, any step-two group can be considered as a subgroup of a class of GM-groups (cf. \cite[\S~8.1]{Li20}).

GM-groups have consummate sub-Riemannian properties, for instance,

\begin{theo}
The following properties are equivalent:
{\em\begin{compactenum}[(i)]
\item $\G$ is of type GM, namely, \[\M= \{ g; \exists \, \theta \in \OA \mbox{ s.t. $\theta$ is a nondegenerate critical point (maximizer) of $\phi(g; \cdot)$ in $\OA$}\} \] is dense in $\G$;
\item $\OM = \left\{\left( x, - \frac{1}{4} \nabla_\theta \langle U(\theta) \, \cot{U(\theta)} \, x, \ x\rangle \right); \ x \in \R^q, \ \theta \in \OA \right\}$ is dense in $\G$;
\item $d(g)^2 = \sup\limits_{\theta \in \OA} \phi(g; \theta)$ for every $g \in \G$;
\item The cut locus of $o$, $\mathrm{Cut}_o$, is exactly the complement of $\M$;
\item The classical cut locus of $o$, $\mathrm{Cut}_o^{\mathrm{CL}} = \OM^c$.
\end{compactenum}}
\end{theo}

We can also characterize GM-group by means of the optimal synthesis from $o$, in other words, the collection of all arclength parametrized geodesics starting from $o$ with their cut times. As a byproduct, the main goal in \cite{BBG12} is achieved from the setting of step-two groups of corank $2$ to all possible step-two groups, via a completely different method. We refer the reader to \cite{Li20} and \cite{LZ21} for further details and more involved results, e.g. the geometric meaning of $\overline{\OM_2}$ as well as the description of shortest geodesic(s) from $o$ to any given $g$ in the setting of GM-groups.

\subsection{Not all step-two groups are of type GM}

Recall that the Gaveau--Brockett problem is completely solved on the free Carnot group of step-two and $3$ generators $N_{3, 2}$ (see \cite{Li20, LZ21}).  The expression of $d(g)^2$ in such case has been provided as explicit as one can possibly hope for, via two highly complicated and non-trivial $C^{\infty}$-diffeomorphism. Thus, the goal of finding an explicit expression for $d^2$ on every step-two group is impractical. Moreover, $N_{3, 2}$ is a non-GM group. Hence the fine characterization of  the squared distance in the framework of GM-groups is no longer valid on $N_{3, 2}$. However, we will be pleasantly surprised to find that the main idea in \cite{Li20}  remains effective to supply a new description of $d^2$, in a uniform way, on any step-two group. For this purpose, more new methods are indispensable. Let us begin with

\subsection{New expressions for the heat kernel} \label{ss4}

Unlike the case of Heisenberg groups or H-type groups, \eqref{2c0} is not enough to study the small-time asymptotics, or sharp estimates for the heat kernel on $N_{3,2}$. Roughly speaking, the reason is that in such situation we have to work on a larger domain  than the initial reference set, which is closely related to the two functions $s \cot{s}$ and $s \csc{s}$ in \eqref{2c0}.

Remark that the main results in \cite{Li20} and \cite{LZ21} only depend on the phase function $\widetilde{\phi}$,  particularly on $s \cot{s}$. In the context of step-two non-GM groups,  the amplitude function  $V$, namely $s \csc{s}$, will be involved in as well and  play an important role. Observe that (cf.  \cite[\S~3.41]{W44}) $\sqrt{\frac{\pi}{2}} \, s^{-\frac{1}{2}} J_{\frac{1}{2}}(s) = \frac{\sin{s}}{s}$, where $J_{c}$ ($c > 0$) denotes the Bessel function. Based on this observation, we will use the properties of Bessel functions to seek more suitable expressions of the heat kernel.

Let $I_{\nu}$ ($\nu > -\frac{1}{2}$) denote the modified Bessel function. It follows from \cite[\S~3.7, \S~15.25, and \S~15.41]{W44} that for each $k \in \N$, the following infinite product representation holds
\begin{align}\label{defQk}
\rQ_k(z) := z^{-\left(k + \frac{1}{2}\right)} \, I_{k + \frac{1}{2}}(z) = \frac{\left( \frac{1}{2} \right)^{k + \frac{1}{2}}}{\Gamma\left(k + \frac{3}{2} \right)} \, \prod_{l = 1}^{+\infty} \left( 1 + \frac{z^2}{\cZ^2_{k,l}}\right),
\end{align}
where $\cZ_{k,1}, \cZ_{k,2}, \ldots$ are the positive zeros of the Bessel function $J_{k + \frac{1}{2}}(\cdot)$, arranged in increasing order, that is, $0 < \cZ_{k,l} < \cZ_{k,l + 1}, \forall \, l \in \N^*$. Remark that from \cite[\S~15.22 and \S~15.4]{W44}, there exists an interlacing property
\begin{align}\label{zerobes}
& 0 < \cZ_{k,1} < \cZ_{k+1, 1} < \cZ_{k,2} < \cZ_{k+1, 2}  < \ldots, \quad \mbox{and} \quad
 \left(\frac{k - 1 }{2} + l \right)  \pi  \le \cZ_{k,l} \le \left(\frac{k + 1 }{2} + l \right)  \pi.
\end{align}
Furthermore, we have the following recursion formulas for modified Bessel functions (cf. \cite[\S~3.71]{W44}):
\begin{align}\label{recbes}
z \, I_{\nu - 1} (z) - z \, I_{\nu + 1}(z) = 2 \, \nu \, I_{\nu}(z) , \qquad
\frac{1}{z} \, \frac{d}{dz} \big\{z^{-\nu} \, I_{\nu}(z) \big\} = z^{- \nu - 1} \, I_{\nu + 1}(z).
\end{align}

Let us introduce for $k \in \N$,
\begin{gather} \label{defRk}
\rR_k(z): = z^2 \, \frac{\rQ_{k + 1}(z)}{\rQ_k(z)} = z \, \frac{\rQ^\prime_k(z)}{\rQ_k(z)}= 2  \, \sum_{l = 1}^{+\infty} \frac{z^2}{\cZ_{k,l}^2 + z^2}, \qquad z \in \C \setminus \{ \pm i \cZ_{k,l}; \, l \in \N^*\}, \\
\label{defPk}
\rPk(X,T):= h^{- \frac{(k + 1)q}{2} - m} \, \int_{\R^m} \left( \det \rQ_k(U(\lambda)) \right)^{-\frac{1}{2}} \, \exp\left\{ -\frac{\langle \rR_k(U(\lambda))  X, \, X \rangle - 2i  T \cdot \lambda }{2h} \right\}  \, d\lambda,
\end{gather}
where $h > 0$  and  $(X,T) \in \R^q \times \R^m$.

Using the well-known result about the Fourier transform of Gaussian functions (cf. \cite[\S~7.6]{H90}), we obtain
\begin{align*}
&e^{-\frac{1}{2h} \, \langle \rR_k(U(\lambda)) \, X, \, X \rangle} =
e^{-\frac{1}{2} \, \left\langle h \, \frac{\rQ_{k + 1}(U(\lambda))}{\rQ_k(U(\lambda))} \, \widetilde{U}(\lambda) \, \frac{X}{h} , \,  \widetilde{U}(\lambda) \, \frac{X}{h} \right\rangle } \\
=& \, \frac{h^{-\frac{q}{2}}}{\Cq} \, \frac{(\det \rQ_k(U(\lambda)))^{\frac{1}{2}}}{(\det \rQ_{k + 1}(U(\lambda)))^{\frac{1}{2}}} \,
\int_{\R^q} e^{-\frac{1}{2h} \, \left\langle \frac{\rQ_k(U(\lambda))}{\rQ_{k + 1}(U(\lambda))} \, s , \, s \right\rangle + \frac{i}{h} \langle \widetilde{U}(\lambda) \, X, \, s \rangle} \, ds.
\end{align*}
Next, taking $\nu = k + 3/2$ in the first equality of \eqref{recbes}, it follows from \eqref{defQk} that
\begin{align}
\rQ_k(z) = (2 k+ 3) \rQ_{k + 1}(z) + z^2 \, \rQ_{k + 2}(z), \ \mbox{namely} \  \frac{\rQ_k(z)}{\rQ_{k + 1}(z)} = (2 k + 3) + z^2 \, \frac{\rQ_{k+ 2}(z)}{\rQ_{k + 1}(z)}.
\end{align}
Combining this with the last integral expression, we deduce the following relation between $\rPk$ and $\rPkk$:
\begin{align}\label{relPk}
\rPk\big(X, T\big) = \frac{1}{\Cq} \, \int_{\R^q} e^{-\frac{2k + 3}{2h} \, |s|^2} \, \rPkk\big(s, T + \langle \UBB \, X, \, s \rangle\big) \, ds.
\end{align}

From our construction and \eqref{2c0}, there exists some constant $\tCqm > 0$ only depending on $q$ and $m$ such that
\begin{align}\label{exphkh}
p_h(x,t) = \tCqm \, e^{-\frac{|x|^2}{4h}} \, \rPo\left(\frac{x}{\sqrt{2}}, t\right).
\end{align}
In the sequel, for $k \in \N^*$, we write an element in $(\R^q)^k$ as $\ss = (\ss_{1} , \ldots, \ss_{k})$ with $\ss_j \in \R^q, 1 \le j \le k$ and define
\begin{align}
\rF_k(x,t,\ss):= \left( \ss_{k}, t + \frac{1}{\sqrt{2}} \, \langle \UBB \, x, \, \ss_1 \rangle +  \sum\limits_{j = 2}^k \langle \UBB \, \ss_{j - 1} , \, \ss_{j} \rangle \right).
\end{align}
Then using \eqref{relPk} repeatedly, we can reformulate the heat kernel into the following more useful expressions:
\begin{align}\label{exphkhk}
p_h(x,t) = \frac{\tCqm}{(\Cq)^k} \, e^{-\frac{|x|^2}{4h}} \, \int_{(\R^q)^k} \exp\left\{- \sum\limits_{j = 1}^k \frac{2j + 1}{2h} \, |\ss_{j}|^2 \right\} \, \rPk\left( \rF_k(x,t,\ss)\right) \, d\ss, \ k = 0,1,2, \ldots.
\end{align}

\subsection{Description of the squared distance on general step-two non-GM groups} \label{ss5}

Intuitively, $\rPk$ shares similar properties with $p_h$. The following result will play a key role below. In particular, we give a direct proof for the positivity of the heat kernel on arbitrary step-two groups by \eqref{exphkh}. It is quite remarkable and somewhat unexpected.

\begin{prop}\label{posPk}
It holds that $\rPk(X,T) > 0$, $\forall \, k \in \N$, $h > 0$, $X \in \R^q$, and $T \in \R^m$.
\end{prop}

To show this proposition, we mainly use \eqref{defQk}, \eqref{defRk},  and the real analyticity w.r.t. the variable $T$ from the Paley--Wiener theorem.

Once we have \eqref{exphkhk} and Proposition \ref{posPk} in hand, we can follow
the same strategy as in \cite{Li20}, and establish the counterparts of/improve those properties (P1)-(P7) above. More precisely, for $k \in \N^*$, we introduce the corresponding reference set and function as follows: $\OAk := \{\tau ; \, \|U(\tau)\| < \cZ_{k,1}\}$,
\begin{gather*}
\rFks(g;\ss,\tau) := |x|^2 +  \sum\limits_{j = 1}^k 2 \, (2j + 1) \, |\ss_j|^2 +\rFk\left(\rF_k(x,t,\ss); \tau \right), \\
 \rFk((X,T); \tau):=  2 \, \langle \rR_k(U(i \tau)) \, X, \, X \rangle + 4 \, T \cdot \tau, \quad \tau \in \OAk,  \quad \mbox{and}\\
 \Mk := \left\{ (X,T) ; \, \exists \, \theta \in \OAk \, \mbox{such that} \, \nabla_\theta \, \rFk((X,T); \theta) = 0, \, \det ( - \He_{\theta} \,  \rFk((X,T); \theta) ) > 0 \right\}.
\end{gather*}

Similarly, $\rFk((X,T); \cdot)$ and $\rFks(g; \ss, \cdot)$  are smooth and concave on $\OAk$. The following counterpart of \eqref{smalltM0} holds: as $h \to 0^+$,
\begin{align}\label{smalltMk}
\rPk(X,T) = \frac{1}{h^{\frac{(k + 1)q + m}{2}}} \frac{(8\pi )^{\frac{m}{2}}}{\left( \det \rQ_k(U(i \theta)) \right)^{\frac{1}{2}}} \, \frac{\exp\{-\frac{\rFk((X,T); \theta)}{4h}\}}{ [\det(-\He_{\theta} \, \rFk((X,T); \theta))]^{\frac{1}{2}}} \, (1 + O(h))
\end{align}
uniformly on every compact set of $\Mk$, where $\theta$ denotes the unique critical point (maximizer) of $\rFk((X,T); \cdot)$ in $\OAk$.

In the following, we set 
\[
\rFzs(g;\tau) := \phi(g;\tau),  \qquad \isrFzs(g) := \srFzs(g) := \sup\limits_{\tau \in \OA} \phi(g;\tau), 
\]
and 
\[
\srFks(g; \ss) := \sup\limits_{\tau \in \OAk} \rFks(g;\ss,\tau), \qquad \isrFks(g) := \inf\limits_{\ss \in (\R^q)^k} \srFks(g;\ss), \qquad k \in \N^*. 
\]
Moreover, we identify $(\R^q)^{k + 1}$ with $(\R^q)^k \times \R^q$ and write an element of $(\R^q)^{k + 1}$ as $(\ss, \ss_{k + 1})$, where $\ss = (\ss_1, \ldots, \ss_k) \in (\R^q)^k$.  The following result can be considered as an improvement of (P4) and (P6) above:

\begin{theo}\label{MTHM}
Let $k \in \N$ and $g = (x,t) \in \G$. We have: \\
(P4') $d(g)^2 \ge \isrFks(g)$; \\
(P6'a) if the function $\rFks(g; \cdot, \cdot)$ attains its minimax, namely  there exists a $(\ss^* , \theta^*) \in (\R^q)^k \times \OAk$  such that
$\isrFks(g) = \srFks(g;\ss^*) = \rFks(g; \ss^*, \theta^*)$,  then $d(g)^2 = \isrFks(g)$ (the meaning of the notations when $k = 0$ herein is obvious); \\
(P6'b) $\isrFks(g) \le \isrFkps(g)$; \\
(P6'c) if $\isrFks(g) = \isrFkps(g)$, then there exists a $(\ss^*,\ss^*_{k + 1}, \theta^*) \in (\R^q)^{k + 1} \times \overline{\OAk}$ such that $\isrFkps(g) = \srFkps(g;\ss^*, \ss^*_{k + 1}) = \rFkps(g; \ss^*, \ss^*_{k + 1}, \theta^*)$, and especially $d(g)^2 = \isrFks(g) = \isrFkps(g)$.
\end{theo}

In particular, if there exists a $k_0$ such that $d(g)^2 = \isrFkzs(g)$, then $d(g)^2 = \isrFks(g), \forall \, k \ge k_0$. Furthermore, from Theorem \ref{MTHM}, we provide an algorithm, which enables us to obtain the exact formula for $d^2$ on $\G$.

\begin{alg}
\textbf{Step (i)} Put in a $g \in \G$ and set $k = 0$;
\textbf{Step (ii)} Compute $\isrFks(g)$;
\textbf{Step (iii)}  If there exists a $(\ss^* , \theta^*) \in (\R^q)^k \times \OAk$  such that
$\isrFks(g) = \srFks(g;\ss^*) = \rFks(g; \ss^*, \theta^*)$, then we obtain $d(g)^2 = \isrFks(g)$ and end this algorithm. If not, return to Step (ii) with $k$ replaced by $k + 1$.
\end{alg}

Finally we must show that the algorithm actually ends in finite steps. Indeed, by the fact that $d^2$ in the setting of step-two groups is locally Lipschitz continuous w.r.t. the usual Euclidean distance (cf.  \cite{R13,R14}), the global reference set
\[
\mathfrak{R} := \overline{\{\theta = 4^{-1} \nabla_t d(x, t)^2; \, g = (x, t) \not\in \mathrm{Cut}_o\}}
\]
is compact. Using the separation theorem of convex sets, we can establish the following:

\begin{theo}\label{MTHM2}
Assume that $k_0 \in \N^*$ such that $\mathfrak{R}  \subset \OAko$. Then $d(g)^2 = \isrFkzs(g)$ for all $g \in \G$.
\end{theo}

\subsection{Sketched proof of Theorem \ref{MTHM}}

Using \eqref{exphkhk},  we can prove (P4') without difficulties by adapting the method in the proof of (P4).  Combining Proposition \ref{posPk} with \eqref{smalltMk} in addition, we can obtain the following upper bound for $d^2$,  which implies immediately (P6'a) via a limiting argument (in some sense, a counterpart of  $\overline{\M} = \overline{\widetilde{\M}}$):

\begin{prop}\label{uppphik}
Let $k \in \N^*$, $g = (x,t) \in \G$, and $\ss^* \in (\R^q)^k$. If $\rF_k(x,t,\ss^*) \in \Mk$, then we have $d(g)^2 \le \srFks(g; \ss^*)$.
\end{prop}

Finally, to show (P6'b) and the first assertion in (P6'c),  we use the results in \cite[\S~36]{R70} and the following:

\begin{lem}\label{mminst}
For $k \in \N$, $g = (x,t) \in \G$, and $\ss^* \in (\R^q)^k$, we have
\begin{align}\label{infrfk+1}
\inf_{\ss_{k + 1} \in \R^q} \rFkps(g;\ss^*,\ss_{k + 1}, \tau) =
\begin{cases}
&\rFks(g;\ss^*,\tau), \qquad  \tau \in \overline{\OAk} \\[2mm]
&-\infty, \qquad \tau \in \OAkp \setminus \overline{\OAk} \\
\end{cases}.
\end{align}
Here we have identified $\rFks(g;\ss^*,\cdot)$ with its closure $\mathrm{cl} \, \rFks(g;\ss^*,\cdot)$ on $\overline{\OAk}$ in the sense of convex analysis (cf. \cite[Section~7]{R70}).
\end{lem}

\subsection{Characterization of normal geodesics from $o$ to any given $g \neq o$ and the cut loci}

We can ameliorate the property (P2).  In particular, for the points belonging to $\mathcal{W}$, using $\nabla_{(\ss,\theta)} \, \rFks (g; \ss, \theta)= 0$ instead of $\nabla_\theta \, \phi (g;\theta) = 0$, we find counterparts of  (P2) and for $k$ large enough, the missing normal geodesics can be found eventually. Furthermore, there exists an one to one correspondence between $\{(g,\zeta,\theta); \, \exp(\zeta, 2 \, \theta) = g \ne o, \theta \in \OAk\}$ and $\{(g,\ss,\theta); \, g \ne o, \nabla_{(\ss,\theta)} \, \rFks(g;\ss,\theta) = 0, \theta \in \OAk\}$.

Recall that the cut locus of the identity $o$, $\mathrm{Cut}_o$, is the set of points where $d^2$ is not smooth. Now, we can characterize it as well. More precisely, setting $k_* := \inf\{k; \, \isrFks(g) = d(g)^2, \forall \, g\} (< + \infty)$, we have:

\begin{theo} \label{T6}
Assume $k > k_*$. Then $g \notin \mathrm{Cut}_o$ if and only if the function $\rFks(g;\cdot,\cdot)$ in $(\R^q)^k \times \OAk$ has a unique critical point $(\ss,\theta)$ in $(\R^q)^k \times \OAk$ such that $\rFks(g;\ss,\theta) = d(g)^2$, and the unique critical point satisfies $\det(\He_{(\ss,\theta)} \, \rFks(g; \ss, \theta)) \ne 0$.
\end{theo}

Similarly, with $k$ as in Theorem \ref{T6}, if both $(\ss,\theta), (\ss^*,\theta^*) \in (\R^q)^k \times \OAk$ (with $\ss \ne \ss^*$)  are simultaneously minimax points and critical points of $\rFks(g;\cdot,\cdot)$, then $g \in \mathrm{Cut}_o^{\mathrm{CL}}$. In addition, it is not hard to provide a counterpart of \cite[Corollary 2.2]{Li20}.

\section{The Riemannian geometry of 2-step groups}

A straightforward modification allows us to obtain most corresponding results of Riemannian geometry related to the full Laplacian. For instance, we can introduce the notation of Riemannian GM-groups and supply their equivalent characterizations via the squared Riemannian distance, the Riemannian cut locus of $o$, as well as the optimal synthesis from $o$.  An interesting and non-trivial fact is that the concept of Riemannian GM-groups is equivalent to that of GM-groups. We can also describe all geodesics from $o$ to any given $g$, and establish the counterparts of  (P7) and Theorem \ref{MTHM}. In addition, the counterpart of Theorem \ref{MTHM2} can still be obtained, even though the scaling property is no longer valid for the Riemannian setting. As an application, we can give an alternative proof for one of two main theorems in \cite{GM03}, namely \cite[Theorem 1 or Theorem 2.9]{GM03} proved therein via Riemannian submersions.

It is worthwhile to point out that usually, the Riemannian cut loci and the sub-Riemannian one do not contain each other. A typical example for such phenomenon is the step-two group of Kolmogorov type on which the reference set is not an open ball.

\section{Other related works}

\subsection{Asymptotic behavior at infinity of the heat kernel}

Compared with \eqref{2c0}, \eqref{exphkhk} (combining with Proposition \ref{posPk}) is more appropriate to reveal the close connection between the geometric properties and the small-time asymptotic behavior for the heat kernel.
Furthermore, by adapting our method, the uniform asymptotic estimate at infinity for the heat kernel as well as sharp bounds for its derivatives can be obtained on a large class of step-two GM-groups, including all GM-M\'etivier groups. Counterparts for the Riemannian heat kernel on GM-M\'etivier groups could be derived as well.

For such topics on non-GM groups, the sharp upper and lower bounds of the heat kernel on $N_{3,2}$ will be provided in a forthcoming paper of us and S.-C. Mao.

\subsection{Riemannian geometry on $N_{3,2}$}

A complete description of the squared Riemannian distance and the Riemannian cut loci will be seen in the thesis of Z.-H. Bao. The proof is closely
related to \cite{Li20}. However, the Riemannian cut loci become much more complicated and more interesting.

\medskip

\section*{Acknowledgement}
\setcounter{equation}{0}
This work is partially supported by NSF of China (Grants  No. 11625102). The authors are indebted to Alano Ancona for the proof of Proposition \ref{NPA}.

\bibliographystyle{abbrv}
\bibliography{LZcc2111}

\mbox{}\\
Hong-Quan Li, Ye Zhang\\
School of Mathematical Sciences/Shanghai Center for Mathematical Sciences  \\
Fudan University \\
220 Handan Road  \\
Shanghai 200433  \\
People's Republic of China \\
E-Mail: hongquan\_li@fudan.edu.cn \\
17110180012@fudan.edu.cn \quad or \quad zhangye0217@126.com \mbox{}\\

\end{document}